\documentclass[12pt]{amsart}
\usepackage{amssymb}
\ifx\mathrm\undefined\let\mathrm\rm\fi
\ifx\mathbf\undefined\let\mathbf\bf\fi
\ifx\mathfrak\undefined\let\mathfrak\frak\fi
\ifx\mathcal\undefined\let\mathcal\cal\fi
\ifx\mathbb\undefined\let\mathbb\Bbb\fi
\ifx\emph\undefined\let\emph\it\fi
\newcommand{\g}{{{\mathfrak g}\,}}

\newcommand{\bor}{{{\mathfrak b}}}
\newcommand{\n}{{{\mathfrak n}}}
\newcommand{\h}{{{\mathfrak h\,}}}

\newcommand{\W}{{{\Bbb W\,}}}
\newcommand{\G}{{{\Bbb G\,}}}
\newcommand{\GH}{{{\Bbb H\,}}}

\newcommand{\Z}{{\mathbb Z}}

\newcommand{\R}{{\mathbb R}}
\newcommand{\C}{{\mathbb C}}

\newcommand{\Ref}[1]{{(\ref{#1})}}
\newcommand{\be}{\begin{displaymath}}
\newcommand{\ee}{\end{displaymath}}
\newcommand{\bea}{\begin{eqnarray*}}
\newcommand{\eea}{\end{eqnarray*}}
\newcommand{\tr}{{\mathrm{tr}}}

\let\dl\delta
\let\Dl\Delta
\let\epe\epsilon \let\eps\varepsilon \let\epsilon\eps

\let\al\alpha

\let\gm\gamma
\let\Gm\Gamma

\let\bt\beta

\let \la \lambda 
\let\om\omega \let\Om\Omega
 \let\phi\varphi

 \let\Si\Sigma

\newcommand{\B}{{\Bbb B}}

\newcommand{\bean}{\begin{eqnarray}}
\newcommand{\eean}{\end{eqnarray}}

\newcommand{\T}{\!\otimes\!\,}

\newcommand{\vs}{\vspace{.5\baselineskip}}
\newtheorem%
{thm}{Theorem}
\newtheorem%
{proposition}[thm]{Proposition}
\newtheorem%
{lemma}[thm]{Lemma}
\newtheorem%
{lemmadef}[thm]{Lemma-Definition}
\newtheorem%
{corollary}[thm]{Corollary}
\newtheorem%
{conjecture}[thm]{Conjecture}
\newcommand{\End}{{\operatorname{End\,}}}

\newcommand{\p}{{\partial}}

\title[ Difference Equations Compatible with KZ Differential Equations]
{  Difference Equations Compatible with Trigonometric KZ Differential Equations}

\author[\smash{V.\;Tarasov and A.\;Varchenko}]
{V.\;Tarasov$^{\,\star,1}$ and
A.\;Varchenko$^{\,\diamond,2}$}
\thanks{$^1$ Supported in part by RFFI grant 99-01-00101}
\thanks{$^2$ Supported in part by NSF grant DMS-9801582}

\begin{document}
\maketitle

\begin{center}
{\it
$^\star$ St. Petersburg Branch of Steklov Mathematical Institute,
Fontanka 27, St. Petersburg 191011, Russia,

vt@pdmi.ras.ru

\medskip

$^\diamond$Department of Mathematics, University of North Carolina
at Chapel Hill,

Chapel Hill, NC 27599-3250, USA,

av@math.unc.edu}
\end{center}
\bigskip
\centerline{February, 2000}
\bigskip
\medskip
\centerline{\sl To the memory of Anatoly Izergin}
\medskip

\begin{abstract}
The trigonometric KZ equations associated with a Lie algebra $\g$
depend on a parameter $\la\in\h$ where $\h\subset\g$ is the Cartan subalgebra.
We suggest a system of dynamical difference equations with respect to $\la$
compatible with the KZ equations. The dynamical equations are constructed
in terms of intertwining operators of $\g$-modules.
\end{abstract}


\thispagestyle{empty}

\section{Introduction}

The trigonometric KZ equations associated with a Lie algebra $\g$
depend on a parameter $\la\in\h$ where $\h\subset\g$ is the Cartan subalgebra.
We suggest a system of dynamical difference equations with respect to $\la$
compatible with the trigonometric KZ differential equations.
 The dynamical equations are constructed
in terms of intertwining operators of  $\g$-modules.

Our dynamical difference equations are a special example of the difference equations introduced
by Cherednik.  In \cite{Ch1, Ch2} Cherednik introduces a notion of an affine R-matrix 
associated with the root system of a Lie algebra and taking values in an algebra $F$ with certain
properties. Given an affine R-matrix, he defines a system of 
equations for an element of the algebra $F$.

In this paper we construct an example of an affine R-matrix and call the corresponding
system of equations the dynamical equations. In our example, $F$ is the algebra of
functions of complex variables $z_1,...,z_n$ and $\la\in\h$ taking values in the tensor product
of $n$ copies of  the universal enveloping algebra of $\g$.
The fact that our dynamical difference equations are compatible with the trigonometric KZ 
differential equations is a remarkable property of our affine R-matrix.

There is a similar construction of dynamical difference equations compatible with
the qKZ difference equations associated with a quantum group. 
The dynamical difference equations in that case are constructed in the same way in
terms of interwining operators of modules over the quantum group.
We will describe this construction in a forthcoming paper.

There is a degeneration of the trigonometric KZ differential equations to
the standard (rational) KZ differential equations. Under this limiting procedure
the dynamical difference equations constructed in this paper turn into
the system of  differential equations compatible with the standard KZ differential equations
and described in \cite{FMTV}. In \cite{FMTV} we proved that the standard hypergeometric
 solutions
of the standard 
KZ equations \cite{SV, V} satisfy also the dynamic differential equations of \cite{FMTV}.

The trigonometric KZ differential equations also have hypergeometric solutions, see \cite{Ch3, EFK}.
We conjecture that the hypergeometric solutions of the trigonometric KZ differential equations
also solve the dynamical difference equations of this paper.

In Section 2 we study relations between intertwining operators of $\g$-modules and the Weyl group $\W$
of $\g$. For any finite dimensional $\g$-module $V$ and $w\in\W$ we construct a rational
function $\B_{w,V} : \C \to \End (V)$. The operators $\B_{w,V}(\la)$ are used later to construct
an affine R-matrix and dynamical equations.

In Section 3 we define the dynamical difference equations for $\g=sl_N$ in terms of operators
$\B_{w,V}(\la)$ directly (without introducing affine R-matrices). For $\g=sl_N$, we 
prove that the dynamical equations are compatible with the trigonometric KZ differential equations.
We give a formula for the determinant of a square matrix solution
of the combined system of KZ and dynamical equations.

In Section 4 we review \cite{Ch1, Ch2} and construct the dynamical difference equations for 
any simple Lie algebra $\g$. We show that the dynamical equations are compatible with
the trigonometric KZ equations if the Lie algebra $\g$ has minuscle weights, i.e. 
is not of type $E_8, F_4, G_2$. We conjecture that the dynamical difference equations and
trigonometric KZ equations are compatible for any simple Lie algebra.

We thank I.Cherednik for valuable discussions and explanation
of his articles \cite{Ch1, Ch2} and P.Etingof
who taught us all about the Weyl group and intertwining operators.

\section{Intertwining Operators}

\subsection{Preliminaries}
Let $\g$ be a complex simple Lie algebra with 
root space decomposition
$\g = \h \oplus(\oplus_{\alpha\in\Si}\g_{\alpha})$
where $\Si\subset\h^*$ is the set of roots.

Fix a system of simple roots $\al_1,...,\al_r$. Let $\Gamma$ be the corresponding
Dynkin diagram, and $\Si_\pm$ --- the set of positive
(negative) roots.
Let $\n_{\pm}=\oplus_{\al\in \Si_{\pm}}\g_\al$. Then
$\g=\n_+\oplus\h\oplus\n_-$.

Let $(\,,\,)$ be
an invariant bilinear form on $\g$.
The form gives rise to a natural identification
$\h\to\h^*$. We use this identification
and make no distinction between $\h$ and $\h^*$.
This identification allows us to define a scalar product on $\h^*$.
We  use the same notation $(\,,\,)$ for the pairing $\h\T \h^*\to\C$.

We use the  notation: $Q=\oplus_{i=1}^r\Z\al_i$ - root lattice;
$Q^+=\oplus_{i=1}^r\Z_{\ge 0}\al_i$;
$Q^\vee=\oplus_{i=1}^r\Z\al_i^\vee$ - dual root lattice,
where $\al^\vee=2\al/(\al,\al)$;
$P=\{\la\in\h\,|\, (\la,\al^\vee_i)\in\Z\}$ - weight lattice;
$P^+=\{\la\in\h\,|\, (\la,\al^\vee_i)\in\Z_{\ge 0}\}$ - cone of dominant integral weights;
$\om_i\in P^+$ - fundamental weights: $(\om_i,\al^\vee_j)=\dl_{ij}$;
$\rho={1\over 2}\sum_{\al\in\Si_+}\al=\sum_{i=1}^r\om_i$;
$P^\vee=\oplus_{i=1}^r\Z\om^\vee_i$ - dual weight lattice, where $\om^\vee_i$
-dual fundamental weights: $(\om^\vee_i,\al_j)=\dl_{ij}$.

Define a partial order on $\h$  putting $\mu<\la$ if $\la-\mu\in Q^+$.

Let $s_i:\h\to\h$ denote 
a simple reflection, defined by $s_i(\la)=\la-(\al_i^\vee,\la)\al_i$;
$\W$ - Weyl group, generated by $s_1,...,s_r$.
The following relations are defining:
\bean\label{rela}
s_i^2=1, \qquad (s_is_j)^m=1 \qquad \text{for}\qquad m=2,3,4,6,
\notag
\eean
where $m=2$ if $\al_i$ and $\al_j$ are not neighboring in $\Gamma$,
otherwise, $m=3,4,6$ if 1,2,3 lines respectively connect $\al_i$ and $\al_j$ in $\Gamma$.
For an element $w\in \W$,  denote $l(w)$ 
the length of the minimal (reduced) presentation of $w$ as a product of generators
$s_1,...,s_r$.

Let $U\g$ be the universal enveloping algebra of $\g$; $U\g^{\T n}$ - tensor product
of $n$ copies of $U\g$; $\Dl^{(n)}:U\g\to U\g^{\T n}$ - the iterated comultiplication
(in particular, $\Dl^{(1)}$ is the identity, $\Dl^{(2)}$ is the comultiplication);
$U\g^{\T n}_0 =\{ x\in U\g^{\T n}\,|\, [\Dl^{(n)}(h),x]=0\,{}\, \text{for any }
h\in\h\}$ - subalgebra of weight zero elements.

For $\al\in\Si$ choose generators $e_\al\in\g_\al$ so that $(e_\al,e_{-\al})=1$.
For any $\al$, the triple 
\bean\label{sl2-al}
H_\al=\al^\vee, \qquad E_\al={2\over (\al,\al)}e_\al,\qquad F_\al=e_{-\al}
\notag
\eean 
forms an $sl_2$-subalgebra in $\g$, $[H_\al,E_\al]=2E_\al,\, [H_\al,F_\al]=-2F_\al,\, 
[E_\al,F_\al]=H_\al$.

A dual fundamental weight $\om_i^\vee$ is called minuscule if $(\om_i^\vee,\al)$ is 0 or 1 for
all $\al\in\Si_+$, i.e. for any positive root $\al=\sum_{i=1}^r m_i \al_i$,
the coefficient $m_i$ is either 0 or 1. For a root system of type $A_r$ all dual
fundamental weights
are minuscule. There is no minuscule dual fundamental weight for $E_8, F_4, G_2$.
For a minuscule dual fundamental weight $\om^\vee_i$, define an element $w_{[i]}=w_0w^i_0 \in \W$
where $w_0$ (respectively, $w_0^i$) is the longest element in $\W$
(respectively, in $\W^i$ generated by all simple reflections $s_j$ preserving
$\om_i^\vee$).

\begin{lemma}\label{wi}
Let $\al$ be a positive root. Then $w_{[i]}(\al)\in\Si_+$ if $(\om_i^\vee,\al)=0$
and $w_{[i]}(\al)\in\Si_-$ if $(\om_i^\vee,\al)=1$.
\end{lemma}

Let $\G$ be the simply connected complex Lie group with Lie algebra $\g$,
$\GH \subset \G$ the Cartan subgroup corresponding to $\h$,
$N(\GH)=\{x\in \G\,|\, x\GH x^{-1}=\GH\}$ the normalizer of $\GH$. Then the Weyl group
is canonically isomorphic to $N(\GH)/\GH$. The isomorphism sends $x$ to
Ad$_x|_\h$.

Let $V$ be a finite dimensional $\g$-module
with weight decomposition $V=\oplus_{\mu\in\h}V[\mu]$.
$\G$ acts on $V$ so that $\GH$ acts trivially on $V[0]$.
Thus the action of $\W$   on $V[0]$ is well defined.
For any $n$, the Weyl group in the same way acts also on $U\g_0^{\T n}$.

\begin{lemma}\label{ef}
For $\al\in \Si$
and $k\in\Z_{\ge 0}$, consider $e_\al^ke_{-\al}^k\in U\g_0$ and $e_{\al}\T e_{-\al}\in
U\g^{\T 2}_0$. Then for any $w\in \W$,
\bean
w(e_\al^ke_{-\al}^k)=e_{w(\al)}^ke_{-w(\al)}^k,\qquad
w(e_\al\T e_{-\al})=e_{w(\al)}\T e_{-w(\al)}.
\notag
\eean
\end{lemma}
\begin{proof}
Let $x\in N(\GH)$ be a lifting of $w$.  Ad$_{x}:\g\to\g$ is an automorphism
of $\g$ preserving the invariant scalar product and sending $\g_\beta$
to $\g_{w(\beta)}$ for all $\beta$. Thus, Ad$_{x}e_{\beta}=c_{x,\beta}e_{w(\beta)}$
for suitable numbers $c_{x,\beta}$ and $c_{x,\al}c_{x,-\al}=1$.
\end{proof}

Let $x_1,...,x_r$ be an orthonormal basis in $\h$, set
\bean
\Om^0={1\over 2}\sum_{i=1}^r x_i\T x_i,
\qquad\Om^+=\Om^0+\sum_{\al\in\Si_+}e_\al\T e_{-\al},\qquad
\Om^-=\Om^0+\sum_{\al\in\Si_+}e_{-\al}\T e_{\al}.
\notag
\eean
Define  the Casimir operator $\Om$ and the trigonometric R-matrix $r(z)$ by
\bean
\Om=\Om^++\Om^- \,,
\qquad 
r(z)={ \Om^+ z+\Om^- \over z-1}\,.
\notag
\eean
For any $x\in U\g$, we have $\Dl(x)\,\Om\,=\,\Om\,\Dl(x)$.
We will use a more symmetric form of the trigonometric R-matrix: $r(z_1/z_2)$.

The Weyl group acts on $r(z_1/z_2), \Om\in U\g^{\T 2}_0$. $\Om$ is Weyl  invariant.
For any $w\in \W$, 
\bean
w(r(z_1/z_2))= {1\over z_1-z_2}\,(\,{z_1+z_2\over 2}\sum_{i=1}^rx_i \T x_i\,+\,
\sum_{\al\in\Si_+} \,(z_1\,e_{w(\al)}\T e_{-w(\al)}\,+\,z_2\,
e_{w(-\al)}\T e_{w(\al)})\,) .
\notag
\eean

\begin{lemma}\label{lemma-2} For a minuscule dual fundamental weight $\om_i^\vee$, 
\bean\label{wr}
z_1^{-(\om_i^\vee)^{(1)}}z_2^{-(\om_i^\vee)^{(2)}}r(z_1/z_2)z_1^{(\om_i^\vee)^{(1)}}z_2^{(\om_i^\vee)^{(2)}}
\,=\, w_{[i]}^{-1}(r(z_1/z_2))\,.
\notag
\eean
\end{lemma}
{\bf Proof.} Using Lemma \ref{wi} 
it is easy to see that both sides of  the equation
are equal to
\bean
{1\over z_1-z_2}\,(\,{z_1+z_2\over 2}\sum_{i=1}^rx_i \T x_i\,+\,
\sum_{\al\in\Si_+,\,(\al,\om^\vee_i)=0} \,( z_1\,e_{\al}\T e_{-\al}\,+\,z_2\,
e_{-\al}\T e_{\al}) +
\notag
\\
 \sum_{\al\in\Si_+,\,(\al,\om^\vee_i)=1} \,(z_1\,e_{-\al}\T e_{\al}\,+\,z_2\,
e_{\al}\T e_{-\al})\,)\,.\qquad \square
\notag
\eean

\subsection{The Trigonometric KZ Equations}

Let $V=V_1\T ...\T V_n$ be a tensor product of $\g$-modules. For $\kappa\in \C$ and
$\la\in \h$, introduce the KZ operators $\nabla_i(\la,\kappa),
\,i=1,...,n,$
acting on functions $u(z_1,...,z_n)$ of $n$ complex variables 
with values in $V$ and defined by
\bean\label{KZ}
\nabla_i(\la,\kappa)
\,=\,
\kappa z_i{\p \over \p z_i}-\sum_{j,\,j\neq i}r(z_i/z_j)^{(i,j)}-\la^{(i)}.
\notag
\eean
Here $r^{(i,j)}$,  $\la^{(i)}$ denote  $r$ acting in 
the $i$-th and $j$-th factors
of the tensor product and $\la$ acting
in the $i$-th factor.


The trigonometric KZ equations are the equations 
\bean\label{KZ-equa}
\nabla_i(\la,\kappa)u(z_1,...,z_n,\la)\,=\,0\,,\qquad i=1,...,n\,,
\eean
 see \cite{EFK}.
The KZ equations are compatible, $[\nabla_i,\nabla_j]=0$.


\subsection{Intertwining Operators, Fusion Matrices, \cite{ES,EV1}}
For $\la\in \h$, let $M_\la$ be the Verma module over $\g$ with highest weight $\la$
and highest weight vector $v_\la$.
We have $\n_+ v_\la=0$, and $\,h v_\la=(h,\la)v_\la$ for
all $h\in\h$. Let $M_\la=\oplus_{\mu\leq\la}M_\la[\mu]$ 
be the  weight decomposition.
The Verma module is irreducible for a generic $\la$.
Define the dual Verma module $M^*_\la$ to be the graded dual space
$\oplus_{\mu\leq \la}M^*_\la[\mu]$  equipped with the  $\g$-action:
$\langle u, a v\rangle=- \langle a u, v \rangle$
for all $a\in\g,\,u\in M_\la,\,v\in M^*_\la$.
Let $v^*_\la$ be 
the lowest weight 
vector of $M^*_\la$ satisfying $\langle v_\la,v^*_\la\rangle =1$.

Let $V$ be a finite dimensional $\g$-module with weight decompostion
$V=\oplus_{\mu\in\h}V[\mu]$. For $\la,\mu\in\h$ consider an intertwining
operator $\Phi\,:\,M_\la\,\to\,M_\mu\T V$.
Define its expectation value by $\langle\Phi\rangle =\langle \Phi (v_\la), v_\mu^*\rangle
\in V[\la-\mu]$.
If $M_\mu$ is irreducible, then the map
Hom$_\g(M_\la,M_\mu\T V)\to V[\la-\mu],\, \Phi\mapsto \langle\Phi\rangle$,
is an isomorphism. Thus for any $v\in V[\la-\mu]$
there exists a unique intertwining operator
$\Phi^v_\la:M_\la\to M_\mu\T V$ such that
$\Phi^v_\la(v_\la)\in v_\la\T v + \oplus_{\nu<\mu}M_\mu[\nu]\T V$.

Let $V,W$ be finite-dimensional $\g$-modules and
$v \in V[\mu],\;w\in W[\nu]$. 
Consider the composition 
\bean
\Phi^{w,v}_{\lambda}:\;M_\lambda
\stackrel{\Phi^v_\lambda}{\longrightarrow} M_{\lambda-\mu}
\otimes V \stackrel{\Phi^w_{\lambda-\mu}}{\longrightarrow}
M_{\lambda-\mu-\nu} \otimes W \otimes V.
\notag
\eean
Then
$\Phi^{w,v}_\lambda \in
\mathrm{Hom}_\g(M_\lambda,M_{\lambda-\mu-\nu}
\otimes W \otimes V)$.  Hence, for a generic $\lambda$
there exists a unique element $u \in \,(V \otimes
W) [\mu+\nu]$ such that
$\Phi^u_\lambda=\Phi^{w,v}_\lambda$. The assignment
$(w, v) \mapsto u$ is bilinear, and defines an $\h$-linear map
$$
J_{WV}(\lambda):\; W \otimes V \to W \otimes V.
$$
The operator $J_{WV}(\lambda)$ is called the fusion matrix of $W$ and $V$.  
The fusion matrix $J_{WV}(\lambda)$ is a
rational function of $\lambda$. $J_{WV}(\lambda)$ is strictly
lower triangular, i.e. $J=1+L$ where $ L(W[\nu] \otimes V[\mu])
\subset \oplus_{\tau<\nu, \,\mu<\sigma} W[\tau]\otimes V[\sigma]$.
In particular, $J_{WV}(\lambda)$ is invertible.

If $V_1,\ldots V_n$ are  $\h$-modules and $F(\lambda):
V_1 \otimes \ldots \otimes V_n \to V_1 \otimes \ldots \otimes V_n$ is a
linear operator depending on $\lambda \in \h$,  then for any
homogeneous $u_1,\ldots , u_n$, $u_i\in V_i[\nu_i]$,
 we  define $ F(\lambda-h^{(i)})(u_1 \otimes
\ldots \otimes u_n)$ to be $F(\lambda-\nu_i) (u_1 \otimes \ldots
\otimes u_n)$.

There is a universal fusion matrix $J(\la)\in 
U\g^{\T 2}_0$ such that  $J_{WV}(\lambda)=J(\lambda)|_{W\otimes V}$ for all $W, V$.
The universal fusion matrix $J(\lambda)$
is the unique solution of the \cite{ABRR} equation
\bean
J(\lambda)\, (1 \otimes (\la+\rho -{1\over 2}\sum_{i=1}^rx_i^2))=
(1 \otimes (\la+\rho -{1\over 2}\sum_{i=1}^rx_i^2) +
\sum_{\alpha \in \Si_+} e_{-\alpha}\otimes e_{\alpha}) J(\lambda).
\notag
\eean
such that $\bigl(J(\la)-1\bigr)\in\bor_-(U\bor_-)\otimes (U\bor_+)\bor_+$
where $\bor_\pm=\h\oplus\n_\pm$.

We transform this equation to a more convenient form. The equation
 can be written as
\bean\label{abrr-1}
J(\lambda)\, (\la+\rho -{1\over 2}\sum_{i=1}^rx_i^2)^{(2)}=
((\la+\rho -{1\over 2}\sum_{i=1}^rx_i^2)^{(2)} -{1\over 2}\sum_{i=1}^rx_i\T x_i
+ \Om^-) J(\la).
\notag
\eean
We make a change of variables: $\la \mapsto  \la -\rho + {1\over 2}(h^{(1)}+h^{(2)})$.
Then the equation  takes the form
\bean\label{abrr-2}
{}\,{}\,{}\,J( \la -\rho + {1\over 2}(h^{(1)}+h^{(2)}))\, 
( \la + {1\over 2}(h^{(1)}+h^{(2)}) -{1\over 2}\sum_{i=1}^rx_i^2)^{(2)}=
\notag
\\
(( \la  + {1\over 2}(h^{(1)}+h^{(2)}) 
-{1\over 2}\sum_{i=1}^rx_i^2)^{(2)} -{1\over 2}\sum_{i=1}^rx_i\T x_i
+ \Om^-) J( \la -\rho + {1\over 2}(h^{(1)}+h^{(2)})).
\notag
\eean
Notice that $(h^{(1)}+h^{(2)})^{(2)}=\sum_{i=1}^r x_i^{(2)}(x_i^{(1)}+x_i^{(2)})$.
Now the equation  takes the form
\bean\label{ABRR}
J( \la -\rho + {1\over 2}(h^{(1)}+h^{(2)}))\,
(\la^{(2)} + \Om^0)\,=\,(\la^{(2)} + \Om^-)\,J( \la -\rho + {1\over 2}(h^{(1)}+h^{(2)})).
\eean

For  $w\in \W$, let $w(J(\la))$ be the image of 
$J(\la)$ under the action of $w$. Let $x\in N(\GH)$ be a lifting of $w$.
Let $W,V$ be finite dimensional $\g$-modules.
Then
\bean\label{wJ}
w(J(\la))|_{W\T V}=x J_{WV}(\la)x^{-1},
\eean
and  RHS does not depend on the choice of $x$.

\subsection{Main Construction, I}\label{main-I}
Introduce a new  action of the Weyl group $\W$ on $\h$
by 
$$
w\cdot \la= w(\la+\rho)-\rho.
$$
Remind facts from \cite{BGG}.

Let $M_\mu, M_\la$ be Verma modules. Two cases are possible:
a) Hom$_\g(M_\mu,M_\la)=0$,
\newline
b) Hom$_\g(M_\mu,M_\la)=\C$ and every nontrivial homomorphism $M_\mu\to M_\la$
is an embedding.

Let $M_\la$ be a Verma module with dominant weight $\la\in P^+$.
Then Hom$_\g(M_\mu,M_\la)=\C$ if and only if there is $w\in \W$
such that $\mu= w\cdot \la$.

Let $w=s_{i_k}\ldots s_{i_1}$ be a reduced presentation.
Set $\al^{1}=\al_{i_1}$ and $\al^{j}=(s_{i_1}\ldots s_{i_{j-1}})(\al_{i_j})$
for $j=2,\ldots,k$. Let $n_j=(\la+\rho,(\al^{j})^\vee)$.
For a dominant $\la\in P^+$, $n_j$ are positive integers.

\begin{lemma}\label{sing-v}
The collection of integers $n_1,\ldots n_k$ and the product
$(e_{-\al_{i_k}})^{n_k}\cdots (e_{-\al_{i_1}})^{n_1}$ do not depend on
the reduced presentation.
\end{lemma}
\begin{proof}
It is known that $\al^{1},\ldots,\al^{k}$ are distinct positive roots and
$\{\al^{1},\ldots,\al^{k}\}=\{\al\in\Si_+\ |\ w(\al)\in\Si_-\}\,$.
Hence, the collection $n_1,\ldots n_k$ does not depend on the reduced
presentation.

The vector $(e_{-\al_{i_k}})^{n_k}\cdots (e_{-\al_{i_1}})^{n_1}v_\la$
is a singular vector in $M_\la$. If $w=s_{i'_k}\ldots s_{i'_1}$ is another
reduced presentation, then the vectors
$(e_{-\al_{i_k}})^{n_k}\ldots (e_{-\al_{i_1}})^{n_1}v_\la$ and
\\
$(e_{-\al_{i'_k}})^{n'_k}\ldots (e_{-\al_{i'_1}})^{n'_1}v_\la$ are
proportional. Since $M_\la$ is a free $\n_-$-module, we have\\
$(e_{-\al_{i'_k}})^{n'_k}\ldots (e_{-\al_{i'_1}})^{n'_1}\,=\,
c\,(e_{-\al_{i_k}})^{n_k}\ldots (e_{-\al_{i_1}})^{n_1}$ in $\n_-$
for a suitable $c\in\C$. $c=1$ since the monomials are equal when projected to
the commutative polynomial algebra generated by $e_{-\al_1},\ldots,e_{-\al_r}$.
\end{proof}

%

Define a singular vector $v_{w\cdot\la}^\la\in M_\la$ by
\bean
v_{w\cdot\la}^\la\,=\,
{(e_{-\al_{i_k}})^{n_k}\over n_1!}
\ldots {(e_{-\al_{i_1}})^{n_1}\over n_k!}\,v_\la\,.
\eean
This vector does not depend on the
reduced presentation
by Lemma \ref{sing-v}.




For all $\la\in P^+$, $w\in \W$, fix an embedding $M_{w\cdot \la}
\hookrightarrow M_\la$ sending $v_{w\cdot \la}$ to
$v_{w\cdot \la}^ \la$. 

Let $V$ be a finite dimensional $\g$-module, $V=\oplus_{\nu\in\h}V[\nu]$
the weight decomposition, $P(V)=\{\nu \in\h\,|\,V[\nu]\neq 0\}$ the set of weights
of $V$. We say that  $\la\in P^+$ is generic with respect to $V$
if
\begin{enumerate}
\item[I.] For any $\nu\in P(V)$ there exist a unique intertwining operator
$\Phi^v_\la:M_\la\to M_{\la-\nu}\T V$ such that $\Phi^v_\la (v_\la)=
v_{\la-\nu}\T v + $  lower order terms. 
\item[II.]
For any $w,w' \in \W,\, w\neq w'$, and any $\nu \in P(V)$,
the vector $w\cdot\la -w'\cdot(\la -\nu)$ does not belong to $P(V)$.
\end{enumerate}
It is clear that all dominant weights lying far inside  the cone of dominant
weights are generic with respect to $V$.

\begin{lemma}\label{a}
Let $\la\in P^+$ be generic with respect to $V$. Let $v\in V[\nu]$.
Consider the intertwining operator $\Phi^v_\la:M_\la\to M_{\la-\nu}\T V$.
For $w\in \W$, consider the singular vector $v_{w\cdot \la}^\la\in M_\la$.
Then there exists a unique vector $A_{w,V}(\la)(v)\in V[w(\nu)]$ such that
\bean\label{A}
\Phi_\la^v (v_{w\cdot\la}^\la)= v_{w\cdot(\la-\nu)}^{ \la-\nu}\T
 A_{w,V}(\la)(v)\,
+\,\text{lower order terms}\,.
\notag
\eean 
\end{lemma}
\begin{proof} $\Phi_\la^v (v_{w\cdot\la}^\la)$ is a singular vector in $M_{\la-\nu}\T V$.
It has to have weight components of
the form $v_{w'\cdot (\la-\nu)}^{\la-\nu}\T u \,$ for suitable $w'\in\W$ and $u\in V$.
Since $\la$ is generic, we have $w=w'$ and $\Phi_\la^v (v_{w\cdot \la}^\la)$ 
is of the required form for a suitable  $A_{w,V}(\la)(v)\in V[w(\nu)]$.
\end{proof}


For  generic $\la\in P^+$,
Lemma \ref{a} defines 
a linear operator $A_{w, V}(\la): V\to V$ such that $A_{w, V}(\la)(V[\nu]))
\subset V[w(\nu)]$ for all $\nu\in P(V)$. It follows from calculations in
Section \ref{main-sl2} that $A_{w, V}(\la)$ is a rational function
of $\la\in\h$.

The following Lemmas are  easy consequences of definitions.

\begin{lemma}\label{Aww}
If $w_1,w_2\in\W$ and $l(w_1w_2)=l(w_1)+l(w_2)$, then
\bea
A_{w_1w_2,V}(\la)\,=\,A_{w_1,V}(w_2\cdot\la)A_{w_2,V}(\la)\,.
\eea
\end{lemma}

\begin{lemma}\label{com-A} 
Let $W, V$ be finite dimensional $\g$-modules. Let $w\in \W$. Then
\bean\label{Com-A}
A_{w,W\T V}(\la)J_{WV}(\la)\,=\,J_{WV}(w\cdot\la)(A_{w,W}(\la-h^{(2)})
\T A_{w,V}(\la))\,.
\notag
\eean
\end{lemma}

Let $x_w\in N(\GH)\subset \G$ be a lifting of $w\in \W$. For a finite dimensional
$\g$-module $V$, define an operator
\bean\label{def-b}
B_{x_w, V}(\la)\,:\,V\,\to\,V\,,\qquad v\,\mapsto \, x_w^{-1}A_{w, V}(\la)v\,.
\notag
\eean
$B_{x_w, V}$ preserves the weight of elements of $V$.

Lemma \ref{com-A} implies
\bean\label{com-xB}
B_{x_w,W\T V}(\la)J_{WV}(\la)\,=\,(x_w^{-1}J_{WV}(w\cdot\la)x_w)\,(B_{x_w,W}(\la-h^{(2)})
\T B_{x_w,V}(\la))\,,
\notag
\eean
cf. \Ref{wJ}.

The operator $B_{x_w, V}$ depends on the choice of $x_w$. If
$x_wg,\, g\in  \GH$, is another lifting of $w$, then
$B_{x_wg, V}\,=\,g^{-1}B_{x_w, V}$.

The operators $B_{x_w,V}(\la)$, $w\in\W$, are defined now for 
generic dominant $\la$
and depend on the choice of liftings $x_w$. 
In the next two Sections we  fix a normalization $B_{w,V}(\la)$
of $B_{x_w,V}(\la)$ so that $B_{w,V}(\la)\,\to \,1 $ as $\la\to\infty$. We  show that
for any $w\in\W$, there is a universal $B_w(\la)\in U\g_0$
such that $B_w(\la)|_V=B_{w,V}(\la)$ for every finite dimensional $\g$-module $V$.
For any $w\in\W$, we  present $B_w(\la)$ as a suitable product of operators $B_{s_i}(\la)$
corresponding to simple reflections.

\subsection{Operators $B_{x_w, V}(\la)$ for $\g=sl_2$}\label{main-sl2}
Consider $sl_2$ with generators $H,E,F$ and relations
$[H,E]=2E,\, [H,F]=-2F,\, [E,F]=H$. Let $\al_1$ be the positive root.
Identifying $\h$ and 
$\h^*$, we have $\al_1=\al_1^\vee=H$,  $\om_1=\om_1^\vee=H/2$,\,
$\W=\{1,s_1\}$. 

Let $\la=l\om_1$, $l\in\Z_{\ge 0}$, be a dominant weight. Then $s_1\cdot \la=
-(l+2)\om_1$. For any dominant weight $\la$, fix an embedding
\bean
M_{s_1\cdot \la}\,\hookrightarrow \,M_\la,\qquad
v_{s_1\cdot\la}\,\mapsto\, v_{s_1\cdot\la}^\la\,=\,{ F^{(\la,\al_1)+1}v_{\la}\over
((\la,\al_1)+1)!}\,
\notag
\eean
as in Section \ref{main-I}.

For $m\in \Z_{\ge 0}$, let $L_m$ be the irreducible $sl_2$ module with highest weight
$m\om_1$. $L_m$ has a basis $v^m_0,..., v^m_m$ such that
$$
Hv^m_k=(m-2k)v^m_k\,,\qquad Fv^m_k=(k+1)v^m_{k+1}\,,\qquad
Ev^m_k=(m-k+1)v^m_{k-1}\,.
$$

For $\g=sl_2$, we have $\G=SL(2,\C)$. Then $\GH\subset\G$ is the subgroup of diagonal
matrices. Fix a lifting ${ x}\in N(\GH)$ of $s_1$,  set $x=(x_{ij})$
where $x_{11}=x_{22}=0$, $x_{12}=-1$, $x_{21}=1$. Then the action of $x$
in $L_m$
is given by $v^m_k\mapsto (-1)^kv^m_{m-k}$ for any $k$.
We have  $x\,=\,\text{exp}(-E)\,\text{exp}(F)\,\text{exp}(-E)$.

For $t\in\C$, introduce
\bean\label{B-sl2}
p(t;\,H,E,F)\,
=
\,\sum_{k=0}^\infty \,F^kE^k\,{1\over k!}\,\prod_{j=0}^{k-1}
{1\over (t-H-j)}\,.
\eean
$p(t;\,H,E,F)$ is an element of $U(sl_2)_0$.

\begin{thm}\label{sl2-B} Let $\la$ be a dominant weight for $sl_2$. Let $L_m,\,x$
be as above. Let $B_{x, L_m}(\la): L_m\to L_m$ be the operator defined in Section
\ref{main-I}. Then for $k=0,...,m$,
\bean\label{Bv}
\\
B_{x, L_m}(\la)v^m_k\,=\,
{((\la,\al^\vee_1)+2)((\la,\al_1^\vee)+3)\cdots ((\la,\al_1^\vee)+k+1)
\over
((\la,\al_1^\vee)-m+k+1)((\la,\al_1^\vee)-m+k+2)\cdots ((\la,\al_1^\vee) - m + 2k)}\,v^m_k\,
\notag
\eean
and 
\bean\label{B=B-un}
p((\la,\al_1^\vee);\,H,E,F)|_{L_m}\,=\,B_{x, L_m}(\la)\,.
\eean
\end{thm}

\begin{corollary}
$B_{x, L_m}(\la)$ is a rational function of $(\la,\al_1^\vee)$.
 $B_{x, L_m}(\la)$ tends to $1$ as $(\la,\al_1^\vee)$ tends to infinity.
\end{corollary}

The Theorem is proved by direct verification. First we calculate explicitly
$\Phi^{v^m_k}_\la \,(v_\la)$, $\Phi^{v^m_k}_\la \,(\,{ F^{(\la,\al_1^\vee)+1}\over
((\la,\al^\vee_1)+1)!}\,v_\la\,)$, and then get
 an expression for $B_{x, L_m}(\la)v^m_k$
as a sum of a hypergeometric type. Using standard formulas from \cite{GR}
we see that $B_{x, L_m}(\la)v^m_k$ is given by \Ref{Bv}. Similarly we check
that $p((\la,\al_1^\vee);\,H,E,F)\,v^m_k$ gives the same result. Thus we
get \Ref{B=B-un}. 
$\square$

Formula \Ref{Bv} becomes 
more symmetric if $\la$ is replaced by $\la-\rho+{1\over 2}\nu$
where $\nu=m\om_1-k\al_1$ is the weight of $v^m_k$, then
\bean\label{product}
p((\la+{1\over 2}\nu,\,\al_1^\vee) - 1;\,H,E,F)
v^m_k\,=\,\prod_{j=0}^{k-1}
{(\la,\al_1^\vee) + {m\over 2} -j
\over (\la,\al_1^\vee) - {m\over 2} +j}\, v^m_k\,.
\eean

\begin{thm}\label{b-sl2}
\bean
p(-t-2;\,-H,F,E)\,\cdot\,p(t;\,H,E,F))\,=
\,{t-H+1\over t+1}\,.
\notag
\eean
\end{thm}
To prove this formula it is enough to check that RHS and LHS give the same result when applied
to $v^m_k\in L_m$, which is done using \Ref{product}. $\square$

Notice that $p(t;\,-H,F,E)=s_1(p(t;\,H,E,F))$.

{\bf Remark.} Let $J(\la)=\sum_ia_i\T b_i$ be the universal fusion matrix of $sl_2$.
Following \cite{EV2}
introduce $S( Q)(\la)\in U(sl_2)_0$ as
$S(Q)(\la)=\sum_i S(a_i)b_i$ where $S(a_i)$ is the antipode of $a_i$. The action of
$S( Q)(\la)$ in $L_m$ was computed in \cite{EV2}. Comparing the result with Theorem
\ref{sl2-B}, one sees that $p((\la,\al_1^\vee);\,H,E,F)$ is equal to
$(S( Q)(\la))^{-1}$ up to a simple change of argument $\la$.

\begin{corollary}\label{AsL}
Let $A_{s_1, L_m}(\la):L_m\to L_m$ be the operator defined in Section \ref{main-I}.
Then $A_{s_1, L_m}(\la)\,=\,x\,p((\la,\al_1^\vee);\,H,E,F)|_{L_m}$.
$A_{s_1, L_m}(\la)$ is a rational function of $(\la,\al_1^\vee)$.
$A_{s_1, L_m}(\la)$ tends to $x$ as $(\la,\al_1^\vee)$ tends to infinity.
\end{corollary}

\subsection{Main Construction, II}\label{main-II}
Return to the situation considered in Section \ref{main-I}.

For any simple root $\al_i$, the triple $H_{\al_i}, E_{\al_i}, F_{\al_i}$
defines an embedding $sl_2 \hookrightarrow \g$ and induces an embedding
$ SL(2,\C) \hookrightarrow \G$. Denote $x_i \in\G$ the image under this embedding
of the element $x\in SL(2,\C)$ defined in Section \ref{main-sl2}.

\begin{lemma}\label{x-in-N}
For $i=1,...,r$, we have $x_i\in N(\GH)$ and Ad$_{x_i}:\g\to\g$ restricted to $\h$
is the simple reflection $s_i:\h\to\h$.
\end{lemma}
\begin{proof} 
Since $x_i\,=\,\exp(-E_{\al_i})$ $\exp(F_{\al_i})$ $\exp(-E_{\al_i})$, we have
that Ad$_{x_i}(H_{\al_i})=-H_{\al_i}$ and Ad$_{x_i}(h)=h$ for any $h\in\h$
orthogonal to $\al_i$. Hence $x_i\in N(\GH)$ and Ad$_{x_i}|_\h=s_i$.
\end{proof}

For $i=1,...,r$ and $\la\in\h$, set 
\bean\label{AB-i}
B_{s_i}(\la)\,=\,p((\la,\al_i^\vee);\,H_{\al_i},E_{\al_i},F_{\al_i})
\notag
\eean
where $p(t;\,H,E,F)$ is defined in \Ref{B-sl2}. Set 
\bean
A_{s_i}(\la)\,=\,x_i\,B_{s_i}(\la)\,.
\notag
\eean
For any $\nu\in P(V)$, we have
$A_{s_i}(\la)(V[\nu])\subset V[s_i(\nu)]$.

Let $V$ be a finite dimensional $\g$-module.
For $w\in\W$, let $w=s_{i_k}...s_{i_1}$ be a reduced presentation.
For a generic dominant $\la\in P^+$, consider the operator $A_{w,V}(\la)
:V\to V$ defined in Section \ref{main-I}. 
\begin{lemma}\label{A=AAA}
\bean
A_{w,V}(\la)\,=\,A_{s_{i_k}}((s_{i_{k-1}}...s_{i_{1}})\cdot\la)|_{V}\,
A_{s_{i_{k-1}}}((s_{i_{k-2}}...s_{i_{1}})\cdot\la)|_{V}...
A_{s_{i_1}}(\la)|_{V}\,.
\notag
\eean
\end{lemma}
\begin{proof}
See Corollary \ref{AsL} and Lemma \ref{Aww}.
\end{proof}
\begin{corollary} The operator $A_{w,V}(\la)$ is a rational function of $\la$.
$A_{w, V}(\la)$ tends to $x_{i_k}...x_{i_1}$ as $\la$ tends to infinity in
a generic direction.
In particular, the product $x_{i_k}...x_{i_1}$ does not depend on the choice
of the reduced presentation.
\end{corollary}
Set $x_w=x_{i_k}...x_{i_1}$. $x_w\in N(\GH)$ is a lifting of $w$. Consider the
operator $B_{x_w,V}(\la):V\to V$ defined in Section \ref{main-I} for this
lifting $x_w$. Denote this operator $B_{w,V}(\la)$.
\begin{corollary}\label{B=BBB}
\bean
&{}&B_{w,V}(\la)\,=
\notag
\\
&{}&(s_{i_{k-1}}...s_{i_{1}})^{-1}(
B_{s_{i_k}}((s_{i_{k-1}}...s_{i_{1}})\cdot\la))|_{V}\,
(s_{i_{k-2}}...s_{i_{1}})^{-1}(B_{s_{i_{k-1}}}((s_{i_{k-2}}...s_{i_{1}})\cdot\la))|_{V}...
B_{s_{i_1}}(\la)|_{V}\,.
\notag
\eean
$B_{w,V}(\la)$ is a rational function of $\la$.
$B_{w, V}(\la)$ tends to $1$ as $\la$ tends to infinity in a generic direction.
\end{corollary}
For any notrivial element $w\in\W$ and $\la\in\h$, define an element $B_w(\la)\in U\g_0$
by
\bean\label{def-B}
&{}&B_{w}(\la)\,=
\notag
\\
&{}&(s_{i_{k-1}}...s_{i_{1}})^{-1}(
B_{s_{i_k}}((s_{i_{k-1}}...s_{i_{1}})\cdot\la))\,
(s_{i_{k-2}}...s_{i_{1}})^{-1}(B_{s_{i_{k-1}}}((s_{i_{k-2}}...s_{i_{1}})\cdot\la))...
B_{s_{i_1}}(\la)\,.
\notag
\eean
Set $B_w(\la)=1$ if $w$ is the identity in $\W$. We have
$B_w(\la)|_V\,=\,B_{w,V}(\la)$, and $B_w(\la)$ does not depend on the choice of
the reduced presentation of $w$.

{\bf Properties of  $B_w(\la)$.}
\begin{enumerate}
\item[I.] 
If $w_1,w_2\in\W$ and $l(w_1w_2)=l(w_1)+l(w_2)$, then
\bean\label{B=BB}
B_{w_1w_2}(\la)\,=\, (w_2)^{-1}(B_{w_1}(w_2\cdot \la))\,B_{w_2}(\la)\,.
\notag
\eean
\item[II.] Let $i=1,...,r$, {}\, $\om\in\h$, and $(\al_i,\om)=0$, then
\bean
B_{s_i}(\la+\om)\,=\,B_{s_i}(\la)\,.
\notag
\eean
\item[III.] For $i=1,...,r$,
\bean
s_i(B_{s_i}(s_i\cdot\la))\,\cdot\,B_{s_i}(\la)\,=\,{(\la,\al_i^\vee)-H_{\al_i}+1
\over (\la,\al_i^\vee)+1}\,.
\notag
\eean
\item[IV.]
Every relation $(s_is_j)^m=1$ for $m=2,3,4,6$ in $\W$ 
is equivalent to a homogeneous relation
$s_is_j...=s_js_i...$. Every such a homogeneous relation  generates 
a relation for $B_{s_i}(\la),
B_{s_j}(\la)$. Namely, for $m=2$, the relation is
\bean\label{ss=ss}
(s_j)^{-1}(B_{s_i}(s_j\cdot\la))\,{}\,B_{s_j}(\la)\,=\,
(s_i)^{-1}(B_{s_j}(s_i\cdot\la))\,{}\,B_{s_i}(\la)\,,
\notag
\eean
for $m=3$, the relation is
\bean\label{ss=ss}
(s_js_i)^{-1}(B_{s_i}((s_js_i)\cdot\la))\,{}\,
(s_i)^{-1}(B_{s_j}(s_i\cdot\la))\,{}\,B_{s_i}(\la)\,=\,&{}&
\notag
\\
(s_is_j)^{-1}(B_{s_j}((s_is_j)\cdot\la))\,{}\,
(s_j)^{-1}(B_{s_i}(s_j\cdot\la))\,{}\,B_{s_j}(\la)\,,
\notag
\eean
and so on.
\item[V.]
\bean\label{com-B}
\Dl (B_{w}(\la))\,J(\la)\,=\,w^{-1}(J(w\cdot\la))\,(B_{w}(\la-h^{(2)})
\T B_{w}(\la))\,.
\notag
\eean
\end{enumerate}

The operators $B_w(\la)$ are closely connected with extremal projectors of
Zhelobenko, see \cite{Zh1, Zh2}.

\subsection{Operators $\B_{w,V}$}\label{Bbb-B}
In order to study interrelations of operators $B_{w,V}(\la)$ with KZ operators
it is convenient to change the argument $\la$.

Let $V$ be a finite dimensional $\g$-module. 
For $w_1,w_2\in\W$ and $\la\in\h$, define 
$w_1(\B_{w_2,V}(\la)):V\to V$ as follows. For any $\nu\in P(V)$ and $v\in V[\nu]$, set
\bean
w_1(\B_{w_2,V}(\la))\,v\,=\,w_1(B_{w_2}(\la-\rho + {1\over 2}\nu))|_{V}\,v\,.
\notag
\eean
In particular,
\bean
\B_{w,V}(\la)v\,=\,B_{w,V}(\la-\rho+{1\over 2}\nu)v\,.
\notag
\eean
$w_1(\B_{w_2,V}(\la))$ is a meromorphic function of $\la$, $w_1(\B_{w_2,V}(\la))$
tends to 1 as $\la$ tends to infinity in a generic direction.

 
{\bf Properties of $\B_{w,V}(\la)$.}
\begin{enumerate}
\item[I.] 
If $w_1,w_2\in\W$ and $l(w_1w_2)=l(w_1)+l(w_2)$, then
\bean\label{B=BB}
\B_{w_1w_2,V}(\la))\,=\, w_2^{-1}(\B_{w_1,V}(w_2( \la)))\,\B_{w_2,V}(\la)\,.
\notag
\eean
\item[II.] If $i=1,...,r$, $w\in\W$,  $v\in V[\nu]$, then
\bean
\B_{s_i,V}(\la)\,v\,
=\,p((\la+{1\over 2}\nu,\,\al_i^\vee)- 1;\,H_{\al_i}, E_{\al_i}, F_{\al_i})\,v
\notag
\eean
and
\bean
w(\B_{s_i,V}(w^{-1}(\la)))\,v\,
=\,p((\la+{1\over 2}\nu,\,w(\al_i^\vee))- 1;\,H_{w(\al_i)}, E_{w(\al_i)}, F_{w(\al_i)})\,v
\notag
\eean
where $p(t;\,H,E,F)$ is defined in \Ref{B-sl2}. 
\end{enumerate}
For $\al\in\Si,\,\la\in\h$, define a linear operator
 $\B^\al_V(\la):V\to V$ by
$$
\B^\al_V(\la)v\,=\,p((\la+{1\over 2}\nu,\al^\vee)-1;\,H_\al,E_\al,F_\al)v
$$
for any $v\in V[\nu]$.

\begin{enumerate}
\item[III.] 
$$
\B^\al_V(\la)\,\B^{-\al}_V(\la)v={(\la-{1\over 2}\nu,\al^\vee)\over
(\la+{1\over 2}\nu,\al^\vee)}v
$$
for any $v\in V[\nu]$.
\item[IV.] Let $\al\in\Sigma$, $\om\in\h$, and $(\al,\om)=0$,
then
\bean
\B_{V}^\al(\la+\om)\,=\,\B_{V}^\al(\la)\,.
\notag
\eean
\item[V.]
Every relation $(s_is_j)^m=1$ for $m=2,3,4,6$ in $\W$ 
is equivalent to a homogeneous relation
$s_is_j...=s_js_i...$. Every such a homogeneous relation  generates 
a relation for $\B_{s_i,V}(\la),
\B_{s_j,V}(\la)$. Namely, for $m=2$, the relation is
\bean\label{ss=ss}
(s_j)^{-1}(\B_{s_i,V}(s_j(\la)))\,{}\,\B_{s_j,V}(\la)\,=\,
(s_i)^{-1}(\B_{s_j,V}(s_i(\la)))\,{}\,\B_{s_i,V}(\la)\,,
\notag
\eean
for $m=3$, the relation is
\bean\label{ss=ss}
(s_js_i)^{-1}(\B_{s_i,V}((s_js_i)(\la)))\,{}\,
(s_i)^{-1}(\B_{s_j,V}(s_i(\la)))\,{}\,\B_{s_i,V}(\la)\,=\,&{}&
\notag
\\
(s_is_j)^{-1}(\B_{s_j,V}((s_is_j)(\la)))\,{}\,
(s_j)^{-1}(\B_{s_i}(s_j(\la)))\,{}\,B_{s_j}(\la)\,,
\notag
\eean
and so on.
\end{enumerate}
These relations can be written in terms of operators $\B^\al_V(\la)$.
\begin{enumerate}
\item[VI.]
For $\al,\beta\in\Si$, denote $\R\langle \al,\bt\rangle $ the subspace $\R\al+\R\bt\subset\h$.
Then
\bean
\B^{\al}_V(\la)\B^{\bt}_V(\la)&=&\B^{\bt}_V(\la)\B^{\al}_V(\la)\,,
\notag
\\
\B^{\al}_V(\la)\B^{\al+\bt}_V(\la)\B^{\bt}_V(\la)&=&
\B^{\bt}_V(\la)\B^{\al+\bt}_V(\la)\B^{\al}_V(\la)\,,
\notag
\\
\B^{\al}_V(\la)\B^{\al+\bt}_V(\la)\B^{\al+2\bt}_V(\la)\B^{\bt}_V(\la)&=&
\B^{\bt}_V(\la)\B^{\al+2\bt}_V(\la)\B^{\al+\bt}_V(\la)\B^{\al}_V(\la)\,,
\notag
\eean
\bean
\B^{\al}_V(\la)\B^{3\al+\bt}_V(\la)\B^{2\al+\bt}_V(\la)\B^{3\al+2\bt}_V(\la)
\B^{\al+\bt}_V(\la)\B^{\bt}_V(\la)=
\notag
\\
\B^{\bt}_V(\la)\B^{\al+\bt}_V(\la)\B^{3\al+2\bt}_V(\la)\B^{2\al+\bt}_V(\la)\,
\B^{3\al+\bt}_V(\la)\B^{\al}_V(\la)\,
\notag
\eean
under the assumption that $\R\langle\al,\bt\rangle =\{\pm\gm\}$ where 
$\gm$ runs over all indices in the
corresponding identity.
\item[VII.]
\bean\label{}
\B_{w,W\T V}(\la))\,=&{}&\,x_w^{-1}(J_{WV}(w(\la)-\rho+{1\over 2}(h^{(1)}+h^{(2)})))x_w\,\cdot
\notag
\\
&{}&(\B_{w,W}(\la-{1\over 2}h^{(2)})
\T \B_{w,V}(\la+{1\over 2}h^{(1)}))\,J(\la-\rho+{1\over 2}(h^{(1)}+h^{(2)}))^{-1}\,
\notag
\eean
\end{enumerate}
\begin{lemma}\label{r-B}
Let $W,V$ be finite dimensional $\g$-modules, $\la\in \h$, $w\in\W$. Then
\bean\label{1}
\Om\,\B_{w,W\T V}(\la)\,=\,
\B_{w,W\T V}(\la) \,\Om
\notag
\eean
and
\bean\label{2}
(w^{-1}(\Om^-)+\la^{(2)})\B_{w,W\T V}(\la)\,=\,
\B_{w,W\T V}(\la) (\Om^-+\la^{(2)})\,.
\notag
\eean
\end{lemma}
\begin{proof}
The first equation  holds since  $\Om$ commutes with the comultiplication.
Now 
\bean
\B_{w,W\T V}(\la) \,(\Om^-+\la^{(2)})\,=
x_w^{-1}(J_{WV}(w(\la)-\rho+{1\over 2}(h^{(1)}+h^{(2)})))x_w\,\cdot
\notag
\\
(\B_{w,W}(\la-{1\over 2}h^{(2)})
\T \B_{w,V}(\la+{1\over 2}h^{(1)}))\,J_{WV}(\la-\rho+{1\over 2}(h^{(1)}+h^{(2)}))^{-1}\,
(\Om^-+\la^{(2)})\,=
\notag
\\
x_w^{-1}(J_{WV}(w(\la)-\rho+{1\over 2}(h^{(1)}+h^{(2)})))x_w\,\cdot
\notag
\\
(\B_{w,W}(\la-{1\over 2}h^{(2)})
\T \B_{w,V}(\la+{1\over 2}h^{(1)}))\,(\Om^0+\la^{(2)})\,J_{WV}(\la-\rho+{1\over 2}(h^{(1)}+h^{(2)}))^{-1}\,=
\notag
\\
x_w^{-1}(J_{WV}(w(\la)-\rho+{1\over 2}(h^{(1)}+h^{(2)})))x_w\,\,(\Om^0+\la^{(2)})\cdot
\notag
\\
(\B_{w,W}(\la-{1\over 2}h^{(2)})
\T \B_{w,V}(\la+{1\over 2}h^{(1)}))\,J_{WV}(\la-\rho+{1\over 2}(h^{(1)}+h^{(2)}))^{-1}\,=
\notag
\\
x_w^{-1}(J_{WV}(w(\la)-\rho+{1\over 2}(h^{(1)}+h^{(2)})))\,(\Om^0+(w(\la))^{(2)})\,x_w\cdot
\notag
\\
(\B_{w,W}(\la-{1\over 2}h^{(2)})
\T \B_{w,V}(\la+{1\over 2}h^{(1)}))\,J_{WV}(\la-\rho+{1\over 2}(h^{(1)}+h^{(2)}))^{-1}\,=
\notag
\\
x_w^{-1}(\Om^-+(w(\la))^{(2)})(J_{WV}(w(\la)-\rho+{1\over 2}(h^{(1)}+h^{(2)})))\,x_w\cdot
\notag
\\
(\B_{w,W}(\la-{1\over 2}h^{(2)})
\T \B_{w,V}(\la+{1\over 2}h^{(1)}))\,J_{WV}(\la-\rho+{1\over 2}(h^{(1)}+h^{(2)}))^{-1}\,=
\notag
\\
(w^{-1}(\Om^-)+\la^{(2)})x_w^{-1}(J_{WV}(w(\la)-\rho+{1\over 2}(h^{(1)}+h^{(2)})))\,x_w\cdot
\notag
\\
(\B_{w,W}(\la-{1\over 2}h^{(2)})
\T \B_{w,V}(\la+{1\over 2}h^{(1)}))\,J_{WV}(\la-\rho+{1\over 2}(h^{(1)}+h^{(2)}))^{-1}\,=
\notag
\\
(w^{-1}(\Om^-)+\la^{(2)})\B_{w,W\T V}(\la)\,.
\notag
\eean
\end{proof}

\section{Difference Equations Compatible with KZ Equations for $\g=sl_{N}$}
\subsection{Statement of Results}
%
Let $e_{i,j}$, $i,j=1,...N$, be the standard generators of the Lie algebra
$gl_N$,
$$
[ e_{i,j}\,, \, e_{k,l}]\,=\,\dl_{j,k}\,e_{i,l}\,-\,\dl_{i,l}\,e_{j,k}\,.
$$
$sl_N$ is the Lie subalgebra of $gl_N$ such that  $sl_n=\n_+\oplus\h\oplus\n_-$ where
$$
\n_+=\oplus_{1\leq i < j\leq N}\C\,e_{i,j}\,,\qquad
\n_-=\oplus_{1\leq j < i\leq N}\C\,e_{i,j}\,, 
$$
and $\h=\{ \la=\sum_{i=1}^N\la_ie_{i,i}\,|\,\la_i\in\C,\,\,\sum_{i=1}^N\la_i=0\}$.

The invariant scalar product is defined by $(e_{i,j}, e_{k,l})=\dl_{i,l}\dl_{j,k}$.
The roots are $e_{i,i}-e_{j,j}$ for $i\neq j$. $\al^\vee=\al$ for any root.
For a root $\al=e_{i,i}-e_{j,j}$, we have $H_\al=e_{i,i}-e_{j,j},\,
E_\al=e_{i,j},\,F_\al=e_{j,i}$.  
The simple roots are $\al_i=e_{i,i}-e_{i+1,i+1}$ for $i=1,...,N-1$.
$\W$ is the symmetric group $S^N$ permutting  coordinates of $\la\in\h$.
The (dual) fundamental weights are $\om_i=\om_i^\vee=\sum_{j=1}^i(1-{i\over N})e_{j,j}
-\sum_{j=i+1}^N{i\over N}e_{j,j}$ for $i=1,...,N-1$. 
All dual fundamental weights are minuscule.
For $i=1,...,N-1$, the permutation $w_{[i]}^{-1}\in S^N$ is
$\left( {}^1_{i+1}\,{}^{2}_{i+2}\,{}^{...}_{...}\,{}^{N-i}_{N}\,{}^{N-i+1}_{1}\,
{}^{...}_{...}\,{}^{N}_{i} \right)$.

For any finite dimensional $sl_N$-module $V$ and $w\in S^N$ consider the operators
$\B_{w,V}(\la)\,:V\to V$. 

Let $V=V_1\T ...\T V_n$ be a tensor product of finite dimensional $sl_N$-modules. For $\kappa\in \C$ and
$\la\in \h$, consider the trigonometric KZ equations with values in $V$,
\bean\label{KZ-sl}
\nabla_j(\la,\kappa)u(z_1,...,z_n,\la)\,=\,0\,,\qquad j=1,...,n\,.
\eean
Here $u(z_1,...,z_n,\la)\in V$ is a function of complex variables $z_1,...,z_n$ and $\la\in\h$.

Introduce {\it the
 dynamical difference equations} on a  
$V$-valued function $u(z_1,...,z_n,\la)$ 
as
\bean\label{dyn-sl}
{}
\\
u(z_1,...,z_n,\la+\kappa \om_i^\vee)\,=\,K_i(z_1,...,z_n,\la)\,
u(z_1,...,z_n,\la)\,, \qquad i=1,...,N-1\,
\notag
\eean
where
$$
K_i(z_1,...,z_n,\la)\,=\,\prod_{k=1}^n z_k^{(\om_i^\vee)^{(k)}}\,\B_{w_{[i]},V}(\la)\,.
$$
The operator $\prod_{k=1}^n z_k^{(\om_i^\vee)^{(k)}}$ is well defined if the argument
of $z_1,...,z_n$ is fixed. The dynamical difference equations are well defined on functions
of $(z,\la)$ where $\la\in\h$ and $z$ belongs to the universal cover of $(\C^*)^n$.
Notice that the KZ equations are well defined for $V$-valued functions of the same variables.

The KZ operators $\nabla_j(\la,\kappa) $ and the operators $K_i(z_1,...,z_n,\la)$ preserve
the weight decomposition of $V$.

\begin{thm}\label{comp-sl}
The dynamical equations \Ref{dyn-sl} together with the KZ equations
\Ref{KZ-sl} form a compatible system of equations.
\end{thm}
\subsection{Proof}
First prove that 
$$
\prod_{k=1}^n z_k^{(\om_i^\vee)^{(k)}}\,\B_{w_{[i]},V}(\la)\,\nabla_j(\la,\kappa) =
\nabla_j(\la+\kappa \om_i^\vee,\kappa)\,\prod_{k=1}^n z_k^{(\om_i^\vee)^{(k)}}\,\B_{w_{[i]},V}(\la)
$$
for all $i$ and $j$. Multiplying both sides from the left by
$\prod_{k=1}^n z_k^{-(\om_i^\vee)^{(k)}}$ and using Lemma \ref{lemma-2}, we reduce the equation to
\bean\label{compa-sl}
\B_{w_{[i]},V}(\la)\,(\,\sum_{k,\,k\neq j} r(z_j/z_k)^{(j,k)}+\la^{(j)}\,)\,=\,
(\,\sum_{k,\,k\neq j}w_{[i]}^{-1}( r(z_j/z_k))^{(j,k)}+\la^{(j)}\,)\,
\B_{w_{[i]},V}(\la)\,.
\notag
\eean
\begin{lemma}\label{nice}
For $j=1,...,n$ and $w\in\W$, we have
\bean\label{compa-sl}
\B_{w,V}(\la)\,(\,\sum_{k,\,k\neq j} r(z_j/z_k)^{(j,k)}+\la^{(j)}\,)\,=\,
(\,\sum_{k,\,k\neq j}w^{-1}( r(z_j/z_k))^{(j,k)}+\la^{(j)}\,)\,
\B_{w,V}(\la)\,.
\notag
\eean
\end{lemma}
\begin{proof}
It is sufficient to check the equation for the residues of both sides
at $z_j=z_k,\,k\neq j$, and for the limit of both sides as $z_j\to\infty$.
The residue equation $[\B_{w,V}(\la), \Om^{(j,k)}]=0$ is true since
the Casimir operator commutes with the comultiplication. The limit equation
\bean\label{compa-sl}
\B_{w,V}(\la)\,(\,\sum_{k,\,k\neq j} (\Om^+)^{(j,k)}+\la^{(j)}\,)\,=\,
(\,\sum_{k,\,k\neq j}w^{-1}(\Om^+)^{(j,k)}+\la^{(j)}\,)\,
\B_{w_{[i]},V}(\la)\,
\notag
\eean
is a corollary of Lemma \ref{r-B}.
\end{proof}

The Theorem is proved for $sl_N,\, N=2$. For $N>2$, it remains to prove  that
\bean\label{compat-dyn}
K_i(z,\la+\kappa\om^\vee_j)\,K_j(z,\la)\,=\,
K_j(z,\la+\kappa\om^\vee_i)\,K_i(z,\la)\,
\eean
for all $i,j$, $0< i<j < N$.
 We prove \Ref{compat-dyn} for $N=3$. For arbitrary $N$
the proof is similar.
 Another proof see in Section \ref{dynamical}.
For $N=3$, $i=1,\, j=2$, equation \Ref{compat-dyn} takes the form
\bean\label{compatib}
{}&{}&\prod_{k=1}^n z_k^{(\om^\vee_1)^{(k)}}
\,\B_{V}^{\al_1+\al_2}(\la+\kappa\om^\vee_2)\,
\,\B_{V}^{\al_1}(\la+\kappa\om^\vee_2)\,
\prod_{k=1}^n z_k^{(\om_2^\vee)^{(k)}}\,
\,\B_{V}^{\al_1+\al_2}(\la)\,
\,\B_{V}^{\al_2}(\la)\,=
\\
&{}&\prod_{k=1}^n z_k^{(\om_2^\vee)^{(k)}}\,
\,\B_{V}^{\al_1+\al_2}(\la+\kappa \om^\vee_1)\,
\,\B_{V}^{\al_2}(\la+\kappa\om^\vee_1)\,
\prod_{k=1}^n z_k^{(\om^\vee_1)^{(k)}}\
\,\B_{V}^{\al_1+\al_2}(\la)\,
\,\B_{V}^{\al_1}(\la)\,.
\notag
\eean
We have $\B_{V}^{\al_1}(\la+\kappa\om^\vee_2)=
\B_{V}^{\al_1}(\la)$
 since $(\om^\vee_2,\al_1)=0$.
We have $[\B_{V}^{\al_1}(\la),
\prod_{k=1}^n z_k^{(\om^\vee_2)^{(k)}}]=0$ since $\B_{V}^{\al_1}(\la)$ is a
power series in $E_{\al_1},\,F_{\al_1}$. Similarly,
$\B_{V}^{\al_2}(\la+\kappa\om^\vee_1)=\B_{V}^{\al_2}(\la)$ and
$[\B_{V}^{\al_2}
(\la),\prod_{k=1}^n z_k^{(\om^\vee_1)^{(k)}}]=0$. Using these remarks and
the relation 
$$
\B_{V}^{\al_2}(\la)
\B_{V}^{\al_1+\al_2}(\la)
\B_{V}^{\al_1}(\la)=
\B_{V}^{\al_1}(\la)
\B_{V}^{\al_1+\al_2}(\la)
\B_{V}^{\al_2}(\la)
$$
we reduce \Ref{compatib} to
\bean
\prod_{k=1}^n z_k^{(\om^\vee_1-\om^\vee_2)^{(k)}}\,
\,\B_{V}^{\al_1+\al_2}(\la+\kappa\om^\vee_2)\,=
\,\B_{V}^{\al_1+\al_2}(\la+\kappa\om^\vee_1)\,
\prod_{k=1}^n z_k^{(\om^\vee_1-\om^\vee_2)^{(k)}}\,.
\notag
\eean
This equation holds since 
$\,\B_{V}^{\al_1+\al_2}(\la+\kappa\om^\vee_2)\,=
\,\B_{V}^{\al_1+\al_2}(\la+\kappa\om^\vee_1)$,
each of these operators is a power
series in  $E_{\al_1+\al_2},\,F_{\al_1+\al_2}$,
 and $(\om^\vee_1-\om^\vee_2, \al_1+\al_2)=0$.

\subsection{An Equivalent Form of Dynamical Equations for $sl_N$}
For $j=1,...,N$, set $\delta_j=\om^\vee_j-\om^\vee_{j-1}$ where
$\om^\vee_0=\om^\vee_N=0$. Then the system of equations
\Ref{dyn-sl} is equivalent to the system
\bean\label{dyn-sl-mod}
u(z_1,...,z_n,\la+\kappa\delta_i)\,=&{}&
\left(\B_V^{e_{i-1,i-1}-e_{i,i}}(\la+\kappa\delta_i)\right)^{-1}...
\left(\B_V^{e_{1,1}-e_{i,i}}(\la+\kappa\delta_i)\right)^{-1}\times
\notag
\\
&{}&\prod_{k=1}^n z_k^{(\delta_i)^{(k)}}\,
\B_V^{e_{i,i}-e_{n,n}}(\la)...\B_V^{e_{i,i}-e_{i+1,i+1}}(\la)
u(z_1,...,z_n,\la)
\notag
\eean
where $i=1,...,N$.  

Notice that the inverse powers 
can be  eliminated using property III in Section \ref{Bbb-B}.

\subsection{Application to Determinants}
Let $\g$ be a simple Lie algebra, $V$ a finite dimensional $\g$-module,
$V[\nu]$ a weight subspace. For a positive root $\al$ fix the $sl_2$ subalgebra
in $\g$ generated by $H_\al, E_\al, F_\al$. Consider $V$ as an
$sl_2$-module. Let $V[\nu]_\al\subset V$ be the $sl_2$-submodule generated by
$V[\nu]$,
$$
V[\nu]_\al=\oplus_{k\in\Z_{\ge 0}}W^\al_k\otimes L_{\nu+k\al}
$$
the decomposition into irreducible $sl_2$-modules. Here
$L_{\nu+k\al}$ is the irreducible module with highest weight
$\nu+k\al$ and $W^\al_k$ the multiplicity space.
Let $d^\al_k=$ dim  $W^\al_k$.
Set
\bean\label{X}
X_{\al,V[\nu]}(\la)=
\prod_{k\in\Z_{\ge 0}} \left(\prod_{j=1}^k
{\Gamma\left(1- {(\la-{1\over 2}(\nu+j\al),\al)\over \kappa}\right)
\over
\Gamma\left(1-{(\la+{1\over 2}(\nu+j\al),\al)\over \kappa}\right)}
\right)^{d^\al_k}\,,
\notag
\eean
cf. formula \Ref{product}. Here $\Gamma$ is the standard gamma function.

Let 
 $V=V_1\T ...\T V_n$ be a tensor product of finite dimensional 
$\g$-modules. Set 
$\Lambda_{k}(\la)=\text{tr}_{V[\nu]}\la^{(k)}$,
$\epe_{k,l}=\text{tr}_{V[\nu]}\Om^{(k,l)}$,
$\gamma_k=\sum_{l,\,l\ne k} \epsilon_{k,l}$.
Set
\bean\label{Det}
D_{V[\nu]}(z_1,...,z_n,\la)=
\prod_{k=1}^nz_k^{{\Lambda_k(\la)\over \kappa} - {\gm_k
\over 2\kappa}}
\,\prod_{1\leq k<l\leq n}(z_k-z_l)^
{\epe_{k,l}\over \kappa}\,
\prod_{\al\in\Si_+}X_{\al,V[\nu]}(\la)\,.
\eean

Let $\g=sl_N$, $V=V_1\T ...\T V_n$ a tensor product 
of finite dimensional $sl_N$-modules.
Fix a basis $v_1,...,v_d$ in a weight subspace $V[\nu]$.  Suppose that
$u_i(z_1,...,z_n,\la)=
\sum_{j=1}^d u_{i,j}v_j$, $i=1,
\ldots,d$, is a set of $V[\nu]$-valued
solutions of the combined system of KZ
 equations \Ref{KZ-sl} and  dynamical equations \Ref{dyn-sl}.

\begin{corollary}\label{determ}
$$
\text{det}\,(u_{i,j})_{1\leq i,j \leq d}\,=\,
C_{V[\nu]}(\la)
\,D_{V[\nu]}(z_1,...,z_n,\la)
$$
where $C_{V[\nu]}(\la)$ is a function of $\la$ (depending also on
$V_1,...,V_n$ and $\nu$) such that
$$
C_{V[\nu]}(\la+\kappa\om)=C_{V[\nu]}(\la)
$$
 for all $\om\in P^\vee$.
\end{corollary}
\begin{proof} The Corollary follows from the following simple Lemma.
\begin{lemma}\label{16}
For $i=1,...,N-1$, the operator $\B_{w_{[i]},V}(\la)$ is the product
in a suitable order of all operators $\B^\al_V(\la)$ with $\al\in\Si_+$
and $(\om^\vee_i,\al)>0$.
\end{lemma}
\end{proof}

Notice that Lemma \ref{16} in particular implies that 
operators $\B_{w_{[i]},V}(\la)$ and the dynamical equations are well
defined in the tensor product of any highest weight $sl_N$-modules.

\section{Dynamical Difference Equations}\label{dynamical}

In this section we introduce dynamical difference equations for arbitrary simple Lie algebra.
The compatibility of the dynamical equations follows from \cite{Ch1}. We prove 
the compatibility of  dynamical and KZ equations.

\subsection{Affine 
Root Systems, \cite{Ch1, Ch2}} Let $\g$ be a simple Lie algebra. The vectors
$\tilde\al = [\al,j]\in \h\times \R$ for $\al\in \Si, j
\in\Z$ form the affine root system $\Si^a$
corresponding to the root system $\Si\subset \h$. We 
view $\Si$ as a subset in $\Si^a$ identifying
$\al\in\h$ with $[\al,0]$. The simple roots of $\Si^a$ are $\al_1,...,\al_r \in \Si$ and
$\al_0=[-\theta, 1]$ where $\theta\in\Si$ is the maximal root. The 
positive roots are
$\Si_+^a=\{[\al,j]\in\Si^a\,|\,
\al\in \Si,\,j>0 \,{}\,\text{or}\, \,\,\al\in\Si_+,\,j=0\}$.
The Dynkin diagram and its affine 
completion with $\{\al_i\}_{0\leq i\leq n}$ as vertices
are denoted
$\Gamma$ and $\Gamma^a$, respectively. The set of the indices of the 
images of $\al_0$ with respect to all authomorphisms
of $\Gm^a$ is denoted $O$ ($O=\{0\}$ for $E_8, F_4, G_2$ ). Let $O^*=\{i\in O\,|\,
i\neq 0\}$. For $i=1,...,r$, the dual
fundamental weight $\om^\vee_i$ is minuscule if and only if
$i\in O^*$. 

Given $\tilde\al=[\al,j]\in\Si^a$ and $\om\in P^\vee$, set
$$
s_{\tilde\al}(\tilde z)=\tilde z - (z,\al^\vee)\tilde\al,
\qquad
t_{\om}(\tilde z)=[z,\xi - (z,\om)]
$$
for $\tilde z=[z,\xi]$.

The affine Weyl group $\W^a$ is the group generated by reflections
 $s_{\tilde\al},\,\tilde\al\in \Si^a_+$.
One defines the length of elements of $\W^a$
taking the simple reflections $s_i=s_{\al_i},\, i=0,...,r$, as generators of $\W^a$.
The group $\W^a$ is the semidirect product $\W \ltimes Q^\vee_t$
of its subgroups $\W=\langle s_\al\,|\,\al\in\Si_+\rangle$ and 
$Q^\vee_t=\{t_\om\,|\, \om\in Q^\vee\}$,
where for $\al\in \Si$ we have $t_{\al^\vee}=s_\al s_{[\al,1]}=s_{[-\al,1]}s_\al$.

Consider the group $P^\vee_t=\{t_\om\,|\, \om\in P^\vee\}$.
The {\it extended affine Weyl group} $\W^b$ is the group of transformations of $\h\times\R$
generated by $\W$ and $P^\vee_t$. $\W^b$ is isomorphic
 to $\W\ltimes P^\vee_t$ with action
$(w,\om)([z,\xi])=[w(z),\xi-(z,\om)]$. 

Notice that for any $w\in\W^b$ and $\tilde\al\in\Si^a$, we have $w(\tilde\al)\in \Si^a$.

The extended affine Weyl group has a remarkable subgroup $\Pi=\{\pi_i\,|\, i\in O\}$,
where $\pi_0\in\Pi$ is the identity element in $\W^b$ and for $i\in O^*$ we have
$\pi_i= t_{\om^\vee_i}w_{[i]}^{-1}$. The group $\Pi$ is isomorphic
 to $P^\vee/Q^\vee$ with
the isomorphism sending $\pi_i$ to the  minuscle weight $\om^\vee_i$.
For $i\in O^*$, the element $w_{[i]}$ preserves the set
 $\{-\theta,\,\al_1,...,\al_r\}$
and $\pi_i(\al_0)=\al_i=w_{[i]}^{-1}(-\theta)$. We have
$$
\W^b=\Pi\ltimes \W^a,
\qquad
\text{where}\qquad
 \pi_is_l\pi^{-1}_i=s_k
\qquad
\text{if}
\qquad
\pi_i(\al_l)=\al_k\,
\qquad
\text{and}\qquad
0\leq k\leq r\,.
$$
We extend the notion of length to $\W^b$. For $i\in O^*,\, w\in \W^a$, we set the length of
$\pi_iw$ to be equal to the length of $w$ in $\W^a$.

\subsection{Affine R-matrices, \cite{Ch1, Ch2}} 
Fix a $\C$-algebra $F$. A set $G=\{G^\al\in F\,|\,\al\in \Si\}$ is called a closed R-matrix
if
\bean
G^{\al}G^{\bt}&=&G^{\bt}G^{\al}\,,
\notag
\\
G^{\al}G^{\al+\bt}G^{\bt}&=&
G^{\bt}G^{\al+\bt}G^{\al}\,,
\notag
\\
G^{\al}G^{\al+\bt}G^{\al+2\bt}G^{\bt}&=&
G^{\bt}G^{\al+2\bt}G^{\al+\bt}G^{\al}\,,
\notag
\\
G^{\al}G^{3\al+\bt}G^{2\al+\bt}G^{3\al+2\bt}
G^{\al+\bt}G^{\bt}&=&
G^{\bt}G^{\al+\bt}G^{3\al+2\bt}G^{2\al+\bt}
G^{3\al+\bt}G^{\al}\,
\notag
\eean
under the assumption that $\al,\bt\in\Si$ and
$\R\langle\al,\bt\rangle =\{\pm\gm\}$ where 
$\gm$ runs over all indices in the
corresponding identity.

A set $G^a=\{\tilde G^{\tilde\al}\in F\,|\,\tilde \al\in \Si^a\}$ is called a closed affine R-matrix
if $\tilde G^{\tilde\al}$ satisfy the same relations where $\al,\bt$ are replaced with
$\tilde\al, \tilde\bt$.

If $G^a$ is an affine R-matrix, for any $w\in\W^b$ define an element
$ \tilde G_w\in F$ as follows. 
Given a reduced presentation $w=\pi_is_{j_l}...s_{j_1}$,
$i\in O$, $0\leq j_1,...,j_l\leq r$, 
set $\tilde G_w=\tilde G^{\tilde \al^l}...\tilde G^{\tilde \al^1}$ where
$\tilde \al^1=\al_{j_1},\, \tilde \al^2=s_{j_1}(\al_{j_2}),\,
\tilde \al^3=s_{j_1}s_{j_2}(\al_{j_3})$,...
The element $\tilde G_w$ does not depend on the reduced presentation of $w$. 
We set $\tilde G_{\text{id}}=1$.

The unordered set $\{\tilde\al^1,...,\tilde\al^l\}$ is denoted $\tilde A (w)$.
There is a useful formula valid for any  (not necessarily minuscule)
dual fundamental weight $\om^\vee_i$, $i=1,...,r$,
\bean\label{useful}
\tilde A(t_{\om^\vee_i})\,=\,\{[\al,j]\,|\,\al\in\Si_+,\,\text{and}\,
(\om^\vee_i,\al)>j\geq 0\}\,,
\eean
Prop. 1.4 \cite{Ch2}.

Introduce the following formal notation: for $w\in\W^b$, $\tilde \al, \tilde\bt\in\Si^a$,
set ${}^w(\tilde G^{\tilde\al})=G^{w(\tilde\al)},\,
{}^w(\tilde G^{\tilde\al}\tilde G^{\tilde\bt})=G^{w(\tilde\al)}G^{w(\tilde\bt)}$,...
Then the elements $\{ \tilde G_w\,|\,w\in\W^b\}$ form a 1-cocycle:
$$
\tilde G_{xy}={}^{y^{-1}}\tilde G_x\,\tilde G_y\,
$$
for all $x,y\in\W^b$ such that $l(xy)=l(x)+l(y)$.

There is a way to construct a closed affine R-matrix if a closed nonaffine R-matrix 
$G=\{ G^{\al}\in F\,|\,\al\in \Si\}$ is given.
Namely, assume that 
the group $P^\vee_t$ acts on the algebra $F$ so that
${}^{t_\om} (G^\al)= G^\al$ whenever $(\om,\al)=0$, $\om\in P^\vee$, $\al\in\Si$.
Then for $\tilde \al=[\al,j]\in\Si^a$, choose $\om\in P^\vee$ so that $(\om,\al)=-j$
and set $\tilde G^{\tilde\al}= {}^{t_\om} (G^\al)$. The set
$G^a=\{\tilde G^{\tilde\al}\in F\,|\,\tilde \al\in \Si^a\}$ is well defined and
forms a closed affine R-matrix called the affine completion of the R-matrix $G$.

 Assume that a closed affine R-matrix $G^a$ is the affine completion
of a closed nonaffine R-matrix $G$. 
Consider the system of equations for an element $\Phi\in F$:
\bean\label{chered}
{}^{t_{-\om^\vee_i}} (\Phi)=\tilde G_{t_{\om^\vee_i}}\Phi\,,
\qquad
i=1,...,r\,,
\eean
where $\om^\vee_1,...,\om^\vee_r$ are the dual fundamental weights. 

\begin{thm}\label{Chered} \cite{Ch1} The system
 of equations \Ref{chered} is compatible,
\bean\label{cher-comp}
{}^{t_{-\om^\vee_i}} (\tilde G_{t_{\om^\vee_j}})\,\,\tilde G_{t_{\om^\vee_i}}\,=
 \,{}^{t_{-\om^\vee_j}} (\tilde G_{t_{\om^\vee_i}})\,\tilde G_{t_{\om^\vee_j}}
\notag
\eean
for $1\leq i<j\leq r$. 
\end{thm}
{\bf Example, \cite{Ch1}.} Let $\al=\al_1,
\,\beta=\al_2,\, a=-\om^\vee_1, \,b=-\om^\vee_2$. Then the system for $A_2$ is
\bean
{}^{t_a}(\Phi)=\tilde G^{\al+\bt}\tilde G^{\al}\Phi,
\qquad
{}^{t_b}(\Phi)=\tilde G^{\al+\bt}\tilde G^{\bt}\Phi.
\notag
\eean
The  system for $B_2$ is
\bean
{}^{t_a}(\Phi)=\tilde G^{\al+2\bt}\tilde G^{\al+\bt}\tilde G^{\al}\Phi,
\qquad
{}^{t_b}(\Phi)=\tilde G^{[\al+2\bt,1]}\tilde G^{\al+\bt}\tilde G^{\al+2\bt} \tilde G^{\bt}\Phi.
\notag
\eean
The system for $G_2$ is 
\bean
{}^{t_a}(\Phi)&=&\tilde G^{[3\al+2\bt,2]}\tilde G^{[3\al+\bt,2]}\tilde G^{[2\al+\bt,1]}
\tilde G^{[3\al+2\bt,1]}\tilde G^{[3\al+\bt,1]}\times
\notag
\\
&{}& \tilde G^{\al+\bt}\tilde G^{3\al+2\bt}\tilde G^{2\al+\bt}
\tilde G^{3\al+\bt}\tilde G^{\al}\Phi,
\notag
\\
{}^{t_b}(\Phi)&=&\tilde G^{[3\al+2\bt,1]} \tilde G^{3\al+\bt}\tilde G^{2\al+\bt}\tilde G^{3\al+2\bt}
\tilde G^{\al+\bt}\tilde G^{\bt}\Phi.
\notag
\eean

\subsection{Affine R-matrix for Dynamical Equations} Fix $\kappa\in\C$ and a
 natural number $n$. Let $F$ be
the algebra of meromorphic functions of $z_1,...,z_n \in \C$ and $\la\in\h$ with values in
$U\g_0^{\T n}$. Define an action of $\W$ on $F$ by
$$
{}^wf(z_1,...,z_n,\la)\,=\, w(f(z_1,...,z_n,w^{-1}( \la)))
$$ 
and an action of $P^\vee_t$ on $F$
by 
$$
{}^{t_\om} f(z_1,...,z_n,\la)\,=\, \prod_{k=1}^n z_k^{\om^{(k)}}f(z_1,...,z_n,\la-\kappa\om)
\prod_{k=1}^n z_k^{-\om^{(k)}}
$$
where $w\in\W$, $\om\in P^\vee$, $f\in F$. 
\begin{lemma} Those actions extend to 
an action of $\W^b=\W\ltimes P^\vee_t$ on $F$, i.e. ${}^{w}({}^{t_\om}f)={}^{t_{w(\om)}}({}^{w}f)$ 
for $w\in\W$, $\om\in P^\vee$, $f\in F$.  $\square$
\end{lemma}

Define a closed nonaffine $F$-valued R-matrix $G_F=\{ G_F^\al\,|\,\al\in\Si\}$ by 
$$
G_F^\al(\la)\,=\,\Delta^{(n)}(p((\la,\al^\vee)-1; H_\al, E_\al, F_\al)).
$$
Properties of operators $\B^\al_V$ described in Section \ref{Bbb-B}
ensure that $G_F$ is a closed R-matrix.  The action of $P^\vee_t$ on $F$ defined above clearly
has the property:
${}^{t_\om} (G_F^\al)= G_F^\al$ whenever $(\om,\al)=0$, $\om\in P^\vee$, $\al\in\Si$.
This allows us to define a closed affine R-matrix
$G_F^a=\{\tilde G_F^{\tilde\al}\in F\,|\,\tilde \al\in \Si^a\}$ as 
the affine completion of the R-matrix $G_F$. Namely,
for $\tilde \al=[\al,j]\in\Si^a$, we choose $\om\in P^\vee$ so that $(\om,\al)=-j$
and set 
$$
\tilde G_F^{[\al,j]}(z_1,...,z_n,\la)\,= \,{}^{t_\om}( G_F^\al)\,=\,
\prod_{k=1}^n z_k^{\om^{(k)}}\,G_F^\al(\la-\kappa\om)\,
\prod_{k=1}^n z_k^{-\om^{(k)}}\,.
$$

Let $V=V_1\T ...\T V_n$ be a tensor product of finite dimensional $\g$-modules.
Let $F_V$ be the algebra of meromorphic functions of $z_1,...,z_n \in \C$ and $\la\in\h$ with values in
\newline
End $(V)$. The closed affine  R-matrix $G^a_F$ induces a closed affine R-matrix
$G^a_V=\{\tilde G_V^{\tilde\al}\}$ where
$$
\tilde G_V^{\tilde\al}(z_1,...,z_n,\la)=\tilde G_F^{\tilde\al}(z_1,...,z_n,
\la+{1\over 2}\sum_{k=1}^n  h^{(k)})|_V\,.
$$
In other words,
$$
\tilde G_V^{[\al,j]}(z_1,...,z_n,\la)\,=\,
\prod_{k=1}^n z_k^{\om^{(k)}}\,\B_V^\al(\la-\kappa\om)\,
\prod_{k=1}^n z_k^{-\om^{(k)}}\,
$$
where $(\om,\al)=-j$
 and the operators $\B^\al_V$ are defined in Section \ref{Bbb-B}.
For any $w\in \W^b$ and $\tilde \al\in\Si^a$, we have
${}^w(\tilde G_V^{\tilde\al})=\tilde G_V^{w(\tilde\al)}$.

Let $\{\tilde G^V_w\in F_V\,|\,w\in\W^b\}$ be the 1-cocycle associated with
the affine R-matrix $G_V^a$. Consider the system 
$$
\prod_{k=1}^n z_k^{-(\om_i^\vee)^{(k)}}\Phi(z_1,...,z_n,\la+\kappa\om^\vee_i)
\prod_{k=1}^n z_k^{(\om^\vee_i)^{(k)}}\,=\,\tilde G^V_{t_{\om^\vee_i}}
(z_1,...,z_n,\la)\Phi(z_1,...,z_n,\la)\,,
$$
$i=1,...,r$, of equations \Ref{chered} associated with the affine R-matrix $G^a_V$.
By Theorem \ref{Chered} this system is compatible.

{\bf Example.} For $\g=sl_N$,  this system of  equations for an element $\Phi\in F_V$ has the form
$$
\prod_{k=1}^n z_k^{-(\om^\vee_i)^{(k)}}\Phi(z_1,...,z_n,\la+\kappa\om^\vee_i)
\prod_{k=1}^n z_k^{(\om^\vee_i)^{(k)}}\,=\,
\B_{w_{[i]},V}(\la)\Phi(z_1,...,z_n,\la)\,,
$$
$i=1,...,N-1$, cf. \Ref{dyn-sl}.

Introduce {\it the dynamical difference
 equations} on a $V$-valued function $u(z_1,...,z_n,\la)$
as
\bean\label{main-en}
{}
\\
\prod_{k=1}^n z_k^{-(\om^\vee_i)^{(k)}}\,u(z_1,...,z_n,\la+\kappa\om^\vee_i)
\,=\, \tilde G^V_{t_{\om^\vee_i}}(z_1,...,z_n,\la)\,u(z_1,...,z_n,\la)\,,
\notag
\eean
$ i=1,...,r$. Notice that the operators 
$\tilde G^V_{t_{\om^\vee_i}}$ 
preserve the weight decomposition of $V$.
Notice also that the operators $\tilde G^V_{t_{\om^\vee_i}}$ are well defined
on the tensor product of any highest weight $\g$-modules
according to formula \Ref{useful}.

An easy corollary of the compatibility of system \Ref{chered} is
\begin{lemma}\label{our-comp}
The dynamical
 difference equations \Ref{main-en} form a compatible system of equations
for a $V$-valued function $u(z_1,...,z_n,\la)$.
\end{lemma}
In particular, for $\g=sl_N$, the Lemma says that the system \Ref{dyn-sl}
is compatible.
\begin{thm}\label{not-thm}
Assume that the Lie algebra $\g$ has a minuscle dual fundamental weight, i.e. $\g$
is not of type $E_8, F_4, G_2$. Then
the dynamical equations \Ref{main-en} together with the KZ equations
\Ref{KZ-equa} form a compatible system of equations. 
\end{thm}
The Theorem is proved in Section \ref{Proof}.

We conjecture that the statement of the
 Theorem holds for any simple Lie algebra.

Let $\g$ be a simple Lie algebra for which the KZ and dynamical equations
are compatible. Let $V=V_1\T ...\T V_n$ be a tensor product 
of finite dimensional $\g$-modules.
Fix a basis $v_1,...,v_d$ in a weight subspace $V[\nu]$.  Suppose that
$u_i(z_1,...,z_n,\la)=
\sum_{j=1}^d u_{i,j}v_j$, $i=1,
\ldots,d$, is a set of $V[\nu]$-valued
solutions of the combined system of KZ
 equations \Ref{KZ-equa} and  dynamical equations \Ref{main-en}.

\begin{corollary}\label{determ-general}
$$
\text{det}\,(u_{i,j})_{1\leq i,j \leq d}\,=\,
C_{V[\nu]}(\la)
\,D_{V[\nu]}(z_1,...,z_n,\la)
$$
where $C_{V[\nu]}(\la)$ is a function of $\la$ (depending also on
$V_1,...,V_n$ and $\nu$) such that
$$
C_{V[\nu]}(\la+\kappa\om)=C_{V[\nu]}(\la)
$$
 for all $\om\in P^\vee$ and $D_{V[\nu]}(z_1,...,z_n,\la)$ is defined in \Ref{Det}.
\end{corollary}
The Corollary follows from formula \Ref{useful}.

\subsection{Proof of Theorem \ref{not-thm}}\label{Proof}
Introduce an action of $\W^b$ on the KZ operators $\nabla_j(\la,\kappa),\,
j=1,...,n$. Namely, for any $w\in\W$, set
$$
{}^w\nabla_j(\la,\kappa)=w(\nabla_j(w^{-1}(\la),\kappa))=
\kappa z_j{\p \over \p z_j}-\sum_{l,\,l\neq j}w(r(z_j/z_l))^{(j,l)}-\la^{(j)}
$$
and for any $\om\in P_t^\vee$ set
\bean
{}^{t_\om}\nabla_j(\la,\kappa)&=&
\prod_{k=1}^n z_k^{\om^{(k)}}
\nabla_j(\la -\kappa\om,\kappa)
\prod_{k=1}^n z_k^{-\om^{(k)}}=
\notag
\\
\kappa z_j{\p \over \p z_j}&-&\prod_{k=1}^n z_k^{\om_i^{(k)}}
\left(\sum_{l,\,l\neq j}r(z_j/z_l)^{(j,l)}\right)
\prod_{k=1}^n z_k^{-\om_i^{(k)}}
-\la^{(j)}\,.
\notag
\eean

The compatibility conditions 
of the dynamical and KZ equations take the form
$$
\tilde G^V_{t_{\om^\vee_i}}(z_1,...,z_n,\la)\,\nabla_j(\la,\kappa)\, =
\,{}^{t_{-\om^\vee_i}}\nabla_j(\la,\kappa)\,
\,\tilde G^V_{t_{\om^\vee_i}}(z_1,...,z_n,\la)
$$
for $i=1,...,r$, $j=1,...,n$.

The compatibility conditions follow from a more general statement.
\begin{thm}\label{Last}
Assume that the Lie algebra $\g$ has a minuscle dual fundamental weight, i.e. $\g$
is not of type $E_8, F_4, G_2$. Then for any $j=1,...,n$ and any $w\in \W^b$ we have
\bean\label{last}
\tilde G^V_{w}(z_1,...,z_n,\la)\,\nabla_j(\la,\kappa)\, =
\,{}^{w^{-1}}\nabla_j(\la,\kappa)\,
\,\tilde G^V_{w}(z_1,...,z_n,\la).
\notag
\eean
\end{thm}
We conjecture that the statement of the
 Theorem holds for any simple Lie algebra.

The Theorem follows from the next four Lemmas.

\begin{lemma} Let $j=1,...,n$. Assume that 
$$
\tilde G^V_{s_l}\nabla_j(\la,\kappa)
=
{}^{s_l}\nabla_j(\la,\kappa)
\tilde G^V_{s_l}\,,
\qquad
{}^{\pi_i}\nabla_j(\la,\kappa)
=
\nabla_j(\la,\kappa)
$$
for $l=0,...,r$ and $i\in O^*$. Then
\bean\label{}
\tilde G^V_{w}(z_1,...,z_n,\la)\,\nabla_j(\la,\kappa)\, =
\,{}^{w^{-1}}\nabla_j(\la,\kappa)\,
\,\tilde G^V_{w}(z_1,...,z_n,\la)
\notag
\eean
for all $w\in\W^b$.
\end{lemma}
\begin{proof}
If $w=\pi_i s_{m_l}...s_{m_1}$ is a reduced presentation, then
\newline
$\tilde G^V_w={}^{s_{m_1}...s_{m_{l-1}}}(\tilde G^V_{s_{m_l}})...
{}^{s_{m_1}}(\tilde G^V_{s_{m_2}})
\tilde G^V_{s_{m_1}}$ and
\bean
\tilde G^V_w\nabla_j(\la,\kappa)=
{}^{s_{m_1}...s_{m_{l-1}}}(\tilde G^V_{s_{m_l}})...
{}^{s_{m_1}}(\tilde G^V_{s_{m_2}})
\tilde G^V_{s_{m_1}}\nabla_j(\la,\kappa)=
\notag
\\
{}^{s_{m_1}...s_{m_{l-1}}}(\tilde G^V_{s_{m_l}})...
{}^{s_{m_1}}(\tilde G^V_{s_{m_2}})
{}^{s_{m_1}}\nabla_j(\la,\kappa)\tilde G^V_{s_{m_1}}=
\notag
\\
{}^{s_{m_1}...s_{m_{l-1}}}(\tilde G^V_{s_{m_l}})...
{}^{s_{m_1}s_{m_2}}\nabla_j(\la,\kappa)
{}^{s_{m_1}}(\tilde G^V_{s_{m_2}})
\tilde G^V_{s_{m_1}}=
\notag
\\
{}^{s_{m_1}s_{m_2}...s_{m_l}}\nabla_j(\la,\kappa)
{}^{s_{m_1}...s_{m_{l-1}}}(\tilde G^V_{s_{m_l}})...
{}^{s_{m_1}}(\tilde G^V_{s_{m_2}})
\tilde G^V_{s_{m_1}}=
\notag
\\
{}^{w^{-1}}\nabla_j(\la,\kappa) \tilde G^V_w\,.
\notag
\eean
\end{proof}

\begin{lemma}\label{very-nice} Let $j=1,...,n$ and $w\in\W$. Then
$$
\tilde G^V_{w}\nabla_j(\la,\kappa)
=
{}^{w^{-1}}\nabla_j(\la,\kappa)
\tilde G^V_{w}\,.
$$
\end{lemma}

\begin{proof} For $w\in\W$ we have $\tilde G^V_w(z_1,...,z_n\la)=\B_{w,V}(\la)$,
and Lemma \ref{very-nice} is equivalent to Lemma \ref{nice}.
\end{proof}

\begin{lemma}\label{very-very-nice} Let $j=1,...,n$ and $i\in O^*$. Then
$$
{}^{\pi_i}\nabla_j(\la,\kappa)
=
\nabla_j(\la,\kappa).
$$
\end{lemma}

\begin{proof} We have $\pi_i=t_{\om^\vee_i}w^{-1}_{[i]}$. Hence
\bean
{}^{\pi_i}\nabla_j(\la,\kappa)={}^{t_{\om^\vee_i}}({}^{w^{-1}_{[i]}}\nabla_j(\la,\kappa))=
{}^{t_{\om^\vee_i}}
(\kappa z_j{\p \over \p z_j}-\sum_{l,\,l\neq j}w^{-1}_{[i]}(r(z_j/z_l))^{(j,l)}-\la^{(j)})=
\notag
\\
\kappa z_j{\p \over \p z_j}-\prod_{k=1}^n z_k^{(\om^\vee_i)^{(k)}}
\left(\sum_{l,\,l\neq j}w^{-1}_{[i]}(r(z_j/z_l))^{(j,l)}\right)
\prod_{k=1}^n z_k^{-(\om_i^\vee)^{(k)}}
-\la^{(j)}\,=\,\nabla_j(\la,\kappa)\,.
\notag
\eean
The last equality follows from Lemma \ref{lemma-2}.
\end{proof}
\begin{lemma}\label{third} Let $j=1,...,n$. 
Assume that the Lie algebra $\g$ has a minuscle dual fundamental weight.
Then
$$
\tilde G^V_{s_0}\nabla_j(\la,\kappa)
=
{}^{s_0}\nabla_j(\la,\kappa)
\tilde G^V_{s_0}\,.
$$
\end{lemma}
\begin{proof} Let $\om^\vee_i$ be a minuscle dual fundamental weight.
We have $s_0=\pi_i^{-1}s_i\pi_i$ and 
$\tilde G^V_{s_0}= {}^{\pi_i^{-1}}(\tilde G^V_{s_i})$
according to the 1-cocycle property.  Now
\bean
{}^{s_0}\nabla_j(\la,\kappa)\tilde G^V_{s_0}=
{}^{\pi^{-1}_is_i\pi_i}\nabla_j(\la,\kappa){}^{\pi_i^{-1}}(\tilde G^V_{s_i})=
{}^{\pi^{-1}_i}({}^{s_i}({}^{\pi_i}\nabla_j(\la,\kappa))\tilde G^V_{s_i})=
\notag
\\
{}^{\pi^{-1}_i}({}^{s_i}(\nabla_j(\la,\kappa))\tilde G^V_{s_i})=
{}^{\pi^{-1}_i}(\tilde G^V_{s_i}\nabla_j(\la,\kappa))=
{}^{\pi^{-1}_i}(\tilde G^V_{s_i}){}^{\pi^{-1}_i}(\nabla_j(\la,\kappa))=
\tilde G^V_{s_0}\nabla_j(\la,\kappa)\,.
\notag
\eean
\end{proof}
Theorems \ref{not-thm} and \ref{Last}  are proved.

\end{document}